\magnification=1200
\baselineskip15pt

\font\bbigbf=cmbx10 scaled \magstep3
\font\bigrm=cmr10 scaled \magstep2

\def\cnl{\centerline}

\def\text{{}}

\def\ff{{\cal F}}
\def\hh{{\cal H}}
\def\ss{{\cal S}}
\def\xx{{\cal X}}
\def\yy{{\cal Y}}
\def\nn{{\cal N}}

\def\Le{\ \le\ }
\def\ignore#1{}
\def\b{\beta}
\def\z{\zeta}
\def\D{\Delta}
\def\sqn{{\sqrt n}}
\def\tL{{\widetilde L}}
\def\tZ{{\widetilde Z}}

\def\Blm{\left|}
\def\Brm{\right|}
\def\Bl{\left(}
\def\Br{\right)}
\def\Blb{\left\{}
\def\Brb{\right\}}

\def\law{{\cal L}}
\def\l{\lambda}

\def\bone{{\bf 1}}

\def\qedbox{\vcenter{\hrule height.5pt\hbox{\vrule width.5pt height8pt
\kern8pt\vrule width.5pt}\hrule height.5pt}}
\def\ep{\hfill$\qedbox$\bsk}

\def\bsk{\bigskip}
\def\re{{\bf R}}
\def\ex{{\bf E}}
\def\pr{{\bf P}}
\def\integ{{\bf Z}}
\def\bone{{\bf 1}}
\def\var{{\rm Var\,}}
\def\sn{\sum_{i=1}^n}
\def\Eq{\ =\ }
\def\ps{\psi}
\def\nin{\noindent}
\def\msk{\medskip}
\def\proof{\nin {\sl Proof.\ }}
\def\be{{\bf e}}
\def\pp{\PP}

\def\BHJ{Barbour, Holst \& Janson}
\def\AGG{Arratia, Goldstein \& Gordon}
\def\BCL{Barbour, Chen \& Loh}

\def\hZ{{\widehat Z}}
\def\tW{{\widetilde W}}
\def\tXi{{\widetilde \Xi}}
\def\tPi{{\widetilde \Pi}}
\def\giv{\,|\,}
\def\uu{{\cal U}}
\def\d{\delta}
\def\Bi{{\rm Bi\,}}
\def\Po{{\rm Po\,}}
\def\Be{{\rm Be\,}}
\def\PP{{\rm PP\,}}
\def\dtv{d_{{\rm TV}}}

\cnl{{\bbigbf  Stein's (magic) method}}\bsk\msk

\cnl{{\bigrm A. D. Barbour and Louis H. Y. Chen}}\bsk
\cnl{{\bigrm Universit\"at Z\"urich and National University of Singapore}}\bsk\bsk 

\cnl{Written in celebration of Charles Stein's 90th birthday in 2010}\bsk\bsk

\nin{\sl 1. Introduction.}\msk
One of the greatest achievements of probability has been its success in approximating the distributions
of arbitrarily complicated random processes in terms of a rather small number of `universal' processes ---
 Brownian motion, the Poisson process, the Ewens sampling formula, Airy processes, stochastic Loewner
 evolution and so on. The standard approach is to consider sequences of processes, indexed by a
 parameter~$n$, and to establish that suitably normalized versions of the processes converge in
 distribution to one of the standard processes as~$n$ tends to infinity. However, for practical purposes,
 it is much more important to know how accurate such an approximation is for a particular process with a
 fixed value of~$n$, and this is a more difficult question to answer. For instance, central limit theorems
 were known already around 1715, and in full generality by 1900, whereas the corresponding approximation
 theorem of Berry and Esseen was only proved in 1941. Stein's method, introduced in 1970, offers a
 general means of solving such problems. By constructing and exploiting a novel characteristic operator
 associated with a random system --- most often, the one used as the approximation --- it turns out to be
 possible to make precise assessments of the approximation error in a wide variety of circumstances.

Stein's original application was in the context of central limit approximation to partial sums of random
 variables having a stationary dependence structure, a problem involving the normal distribution and the
 real line. However, his method has a big advantage over most other techniques, in that it can in
 principle be used for approximation in terms of any distribution on any space, including random
 variables on the real line, processes on a space of sequences, functions or measures, and combinatorial
 structures on discrete spaces. A further big advantage over its competitors is that strong independence
 requirements are not needed to make the method work (though they may of course simplify many arguments
 and the form of the bounds that can be attained). As a result of this considerable freedom, its uses
 have proliferated, with approximations not only to the normal distribution, but also to the Poisson
 distribution, to multivariate normal distributions, to diffusions, to Poisson processes, to the Ewens
 sampling formula, to the Wigner semi-circle law, and more. 

The method continues to produce new developments in a wide variety of settings.
The last five years alone have seen great progress in problems concerning large deviations and 
concentration of measure inequalities; in the application of
 Stein's method to problems having an essentially algebraic component; 
in proving bounds for normal and gamma approximations to the distributions of functionals of 
infinite--dimensional Gaussian fields, using a combination of Stein's method and Malliavin calculus;
and in a range of problems involving random geometrical graphs. It is a tribute to the importance
of Stein's original and amazing idea that it continues to inspire vigorous research 40 years after
its inception. In this article, we can do no more than scratch
the surface of a wide-ranging and still mysterious topic.

\bsk\msk

\nin{\sl 2. Normal approximation}\msk
Stein originally introduced his method in a course of lectures at Stanford.
Dissatisfied with the then available proofs of the combinatorial central limit
theorem, he devised a new one for his class, which entirely dispensed with the
use of Fourier analysis. His paper (Stein, 1972) in the VI'th Berkeley Symposium
contained the first published version of his method, giving Berry--Esseen bounds 
for the accuracy of the normal approximation to the distribution of a sum
of a stationary sequence of random variables; already, the
flexibility of the method is in evidence, since the dependence setting is quite
different from that of the combinatorial central limit theorem, with which
he began.  One way of arriving at his new approach to normal approximation is as follows.

Let $Z$ be a random variable on~$\re$ with differentiable probability density~$p$.
Our main emphasis here is for~$p$ the standard normal density, but this argument
works more generally. 
Now let~$\ff$ be the set of differentiable real functions~$f$ such that $f(x)p(x) \to 0$ as 
$|x| \to \infty$.  Then, clearly,
$$\eqalign{
  0 &\Eq \int_{-\infty}^\infty {d\over dx}\{f(x)p(x)\} \, dx \cr
    &\Eq \int_{-\infty}^\infty \{f'(x) + \ps(x)f(x)\}p(x)\,dx,} \eqno(2.1)
$$
with $\ps(x) := p'(x)/p(x)$: for~$p$ the standard normal density, $\ps(x) = -x$.   
Hence any function~$h$ of the form $f'+\ps f$ with~$f\in\ff$ automatically has $\ex h(Z) = 0$.  
Conversely, given any continuous~$h$ with $\ex|h(Z)| < \infty$,
the function $h - \ex h(Z)$ can be written in the form $f'+\ps f$ with~$f\in\ff$.
To do so is simple.  First observe that, if~$g = f'+ \ps f$, then 
$$\eqalign{
   \int_{-\infty}^X g(x)p(x) \, dx &\Eq \int_{-\infty}^X {d\over dx}\{f(x)p(x)\} \, dx 
            \Eq f(X)p(X),  }       \eqno(2.2)
$$
and that this argument can be reversed: $f$ so defined, for an arbitrary~$g$ for which 
$\ex |g(Z)| < \infty$, is such that $g = f' + \ps f$.
Hence we can take~$f=f_h$ to be defined by
$$\eqalign{
   f_h(X)p(X) 
            &\Eq \int_{-\infty}^X \{h(x) - \ex h(Z)\}p(x) \, dx \cr
            &\Eq - \int_X^{\infty} \{h(x) - \ex h(Z)\}p(x) \, dx,} \eqno(2.3)
$$
noting that it is then directly checked that $f_h\in\ff$.  This allows us to write
$$
   \ex h(W) - \ex h(Z) \Eq \ex\{f_h'(W) + \ps(W) f_h(W)\}, \eqno(2.4)
$$
for any random variable~$W$ for which the expectations exist; in particular, for standard
normal approximation,
$$
   \ex h(W) - \ex h(Z) \Eq \ex\{f_h'(W) - Wf_h(W)\}. \eqno(2.5)
$$  

Taking the supremum
of the left hand side in~(2.4) over test functions~$h$ in some suitable class~$\hh$ 
gives a (very
concrete) measure of the distance between the distributions of $W$ and~$Z$.
This distance can in turn be computed by taking the supremum of the right
hand side of~(2.4) over $h\in\hh$.  The intuition is then that, since the right hand side 
is exactly zero for {\it all\/} $f\in\ff$ if $\law(W) = \law(Z)$, it should
`automatically' be close to zero if $\law(W) \approx \law(Z)$, again more or
less irrespective of~$f$, and hence the supremum of the right hand side of~(2.4)
for $h\in\hh$ may indeed be shown to be small, as it were by right.  The
remarkable power of the method derives from the fact that, in many many 
circumstances, expressing the difference $\ex h(W) - \ex h(Z)$ in the
entirely equivalent form $\ex\{f_h'(W) + \ps(W) f_h(W)\}$ makes it much
easier to bound.

A typical setting, reasonable for standard normal approximation, is one in which~$W$ is
a sum $\sn X_i$ of many individually small and only weakly dependent random
variables~$X_i$, with zero means and with $\var W = 1$. Since, for standard
normal approximation,
$\ps(W) = - \sn X_i$, the expression $\ex \{\ps(W) f_h(W)\}$ can be broken
up into a sum of terms $\ex\{-X_i f_h(W)\}$, in which $X_i$ has only limited
influence on the whole sum~$W$. This can then be exploited, together
with a Taylor expansion, to show that $\ex\{f_h'(W) + \ps(W) f_h(W)\}$ is
small.  For instance, for {\it independent\/} summands~$X_i$, one can
write
$$\eqalign{
   \ex\{W &f_h(W)\} \Eq \sn \Bl\ex\{X_i f_h(W-X_i)\} + \ex\{ X_i[f_h(W) - f_h(W-X_i)]\} \Br\cr
       &\Eq \sn\Bl\ex\{X_i^2 f'_h(W-X_i)\} + \ex\{ X_i[f_h(W) - f_h(W-X_i) - X_if'_h(W-X_i)]\}\Br,
}$$
using the independence of $W-X_i$ and~$X_i$ with $\ex X_i = 0$, and compare the result with
$$
   \ex f'_h(W) \Eq \sn \ex\{X_i^2\} \ex f'_h(W),
$$
because $\ex W^2 = 1$.  This immediately gives
$$
   |\ex h(W) - \ex h(Z)| \Eq |\ex\{f'_h(W) - Wf(W)\}|
     \Le {3\over2}\sn \ex|X_i|^3 \|f''_h\|_\infty,  \eqno(2.6)
$$
yielding an explicit Lyapounov bound for classes of test functions~$\hh$ for which 
the supremum 
$\sup_{h\in\hh} \|f''_h\|_\infty < \infty$.  This elementary argument can of
course be improved in many ways, to deal with less restrictive classes of test 
functions: see Section~6.  What is more important to notice is that independence
can immediately be replaced by some form of local dependence in the argument,
without any great change in spirit, by using the dissection
$$
  X_i f_h(W) \Eq X_i f_h(W-Y_i) +  X_i[f_h(W) - f_h(W-Y_i)],
$$
where now $X_i$ and $W-Y_i$ are (almost) independent and $Y_i$ is still small
enough not to have too big an influence on~$W$.  This illustrates that the
method can be much more robust with respect to dependence than the classical
transform methods.

The discussion above is intentionally kept as simple as possible, and this to
some extent disguises the power of the method.  Even for normal approximation, 
Stein's method has had many startling successes.  One of the earliest was
Bolthausen's~(1984) Berry--Esseen Lyapounov bound for the error in the 
combinatorial central theorem, which had been the subject of much research
in the previous decades.  G\"otze's~(1991)
multivariate normal approximation theorem is another: using a proof based
on Stein's method, he was able to establish 
an error bound of order $O(n^{-1/2})$ for the probabilities of an extremely large
collection of sets, including all convex sets, for sums of independent random
vectors.  Rinott \& Rotar~(1996) developed ideas originating in Stein's proof of (5.3)
and in G\"otze~(1991), in
proving an explicit and very effective error bound, typically of order $n^{-1/2}\log n$,
 for the normal approximation of 
sums of bounded, dependent random vectors $X_j$ in $\re^d$. The dependence structure envisaged 
there
is one in which only relatively few other ``neighbouring'' random vectors are significantly dependent 
on any given $X_j$, with the meaning of neighbouring to be chosen to suit:
see also Chen \& Shao~(2004).

\bsk\msk
\nin{\sl 3. Stein's method in a nutshell.}\msk
 Stein's method can be simply described as follows. Let $W$ and~$Z$ be
 random elements taking values in a space~$\ss$, and let $\xx$ and~$\yy$ be classes of 
bounded real--valued functions defined on~$\ss$. In approximating the 
distribution~$\law(W)$ of~$W$ by the distribution 
$\law(Z)$ of~$Z$, we write 
$\ex h(W) - \ex h(Z) = \ex\{Lf_h(W)\}$ for $h\in \yy$, where $L$ is a linear operator 
from $\xx$ into~$\yy$, and~$f_h$ a bounded solution of 
the equation $Lf = h - \ex h(Z)$. The error $\ex \{Lf_h(W)\}$ can then be bounded by studying the 
solution~$f_h$ and by exploiting the probabilistic properties of~$W$.

Usually, we take as test functions a class $\hh\subset\yy$, large enough to be separating; 
that is, large enough so that, if $\ex h(W) = \ex h(Z)$ for all $h\in\hh$, then 
$\law(W) = \law(Z)$. $\xx$ is then assumed to contain all $f_h$ for   $h\in\hh$.  In this case, 
$\ex \{Lf(W)\} = 0$ for all $f\in\xx$ if
 and only if $\law(W) = \law(Z)$.  Such an $L$ characterizes~$\law(Z)$; the equation 
$Lf = h - \ex h(Z)$ is called a 
Stein equation for~$\law(Z)$, and $L$ a Stein operator for~$\law(Z)$.
In the particular case of normal approximation, where $W$ and~$Z$ are real-valued random variables 
and $\law(Z) = \nn(0, 1)$, 
the operator~$L$ used by Stein~(1972) is given by $Lf(w) = f'(w) - wf(w)$, as noted
above.

\bsk\msk
\nin{\sl 4. Stein identities.}\msk
The discussion so far has been rather vague as to why it might be easier
to show that $\ex \{Lf_h(W)\}$ is small than it is to show the same for $\ex h(W) - \ex h(Z)$
directly.
For the normal distribution, the motivation came by way of an analytic plausibility
argument, and it was demonstrated by example that it was true for
the chosen Stein operator~$L$.  However,
it is implied in Stein~(1992) that one can directly exploit the probabilistic properties 
of~$W$, using auxiliary randomization, in order to find a linear operator~$L$ for 
which the error $\ex \{Lf_h(W)\}$ is
manageable, and indeed thereby deducing the distribution~$\law(Z)$ that should be used
to approximate~$\law(W)$.   The idea is to begin by finding an operator~$\tL$  
such that $\ex \{\tL f(W)\} = 0$ for all $f\in\xx$;  that is, to begin by finding 
a Stein identity for~$\law(W)$.  Once this has been done, one can look for an operator~$L$
that characterizes a better known distribution~$\law(Z)$, with the property that $Lf$ is close
to~$\tL f$ for all $f\in\xx$. The error $\ex \{Lf_h(W)\}$ in the approximation is then 
equivalently expressed as $\ex \{Lf_h(W) - \tL f_h(W)\}$, because $\ex \{\tL f(W)\} = 0$
for all~$f$,
and this difference is small because $Lf_h$ is close to~$\tL f_h$.

This approach is most simply illustrated when the random variable~$W$ is known to have the
equilibrium distribution of some Markov process~$\tZ$ with generator~$\tL$.  In that
case, for suitable functions~$f$, $\ex \{\tL f(W)\} = 0$ by Dynkin's formula, and if~$\tL$ 
is close to the generator~$L$ of another Markov process whose limit distribution~$\law(Z)$
is well known, approximations can be deduced.  More frequently, $W$ has the equilibrium
distribution of some function~$g(\tZ)$, so that $\tL f(W)$ is replaced by
$\tL(f\circ g)(\tZ)$, which need not be a function of~$g(\tZ)$ alone; however, it may still 
be possible to find an operator~$L$ on~$\xx$ such that $Lf(g(\tZ))$ is close to
$\tL(f\circ g)(\tZ)$, which is all that is needed to  make the method work.
Indeed, having guessed such an~$L$, one can then consider the expressions $\ex\{Lf_h(W)\}$ 
directly, with a well founded hope that they can be shown to be small.  This idea has
become known as the generator method (Barbour 1988).  One of its greatest successes has been the 
multivariate normal approximation theorem of G\"otze~(1991), referred to above.  
The use of the generator method in Poisson process approximation is discussed 
in Section~7,
and its application to diffusion approximation can be found in Barbour~(1990).

As a simpler example, suppose that $W = \sn X_i$, with the $X_1,\ldots,X_n$  independent and 
having zero means, and with $\ex W^2=1$.  Define a Markov chain~$\tZ$ as follows.  Let
$I_j$, $j\ge1$, be independent (also of the $X_i$'s) and uniformly distributed 
on $\{1,2,\ldots,n\}$, and let $\tZ(0) = (X_1,\ldots,X_n)$.  Then, if $\tZ(j-1) = z$, 
define $\tZ(j) = z + \be_{I_j}(X_{I_j}^{(j)} - z_{I_j})$, where the $X_{i^{(j)}}$, $j\ge1$,
are independent copies of~$X_i$, $1\le i\le n$, and all are independent of everything else:
$\be_i$ denotes the $i$-th unit vector.
Clearly, $\tZ$ is a stationary Markov chain, and its distribution at any time is that of 
independent random variables with the distributions of the~$X_i$.
The generator is given by
$$
   \tL f(z) \Eq {1\over n}\sn \ex\{f(z + \be_i(X_i-z_i)) - f(z)\}. \eqno(4.1)
$$

We are actually interested in $W = g(\tZ(0))$, where $g(z) = \sn z_i$. For functions
$f\circ g$ for this~$g$, we have
$$\eqalign{
   \tL(f\circ g)(z) &\Eq {1\over n}\sn \ex\{f(g(z + \be_i(X_i-z_i))) - f(g(z))\} \cr
     &\Eq {1\over n}\sn \ex\Blb (X_i-z_i)f'(g(z)) + {1\over2}(X_i-z_i)^2 f''(g(z)) \Brb
       + {1\over n}E(f, g,z), }     \eqno(4.2)
$$
where
$$\eqalign{
    &E(f, g,z) \cr
    &\Eq \sn \ex\Blb f(g(z + \be_i(X_i-z_i))) - f(g(z))  
        - (X_i-z_i)f'(g(z)) - {1\over2}(X_i-z_i)^2 f''(g(z)) \Brb. } \eqno(4.3)
$$
Using the fact that $\ex X_i = 0$ for each~$i$ and that $\ex W^2 = 1$, and
absorbing the factor~$1/n$ (which plays no real part) into the definition of~$\tL$,
this gives
$$\eqalign{
   \tL(f\circ g)(z) &\Eq  -g(z)f'(g(z)) + {1\over2}\Bl 1 + \sn z_i^2\Br f''(g(z))
     + E(f,g,z).   }  \eqno(4.4)
$$
This shows that $\tL(f\circ g)(z)$ is close to $Lf(g(z))$, with the Stein operator
given by
$$
     Lf(w) \Eq -wf'(w) + f''(w), 
$$
provided that $\Bl \sn z_i^2 - 1\Br$ and $E(f,g,z)$ are small.  
Note that the Stein operator is the same as that derived earlier for the
standard normal distribution,
except that the derivative~$f'$ now plays the part of~$f$. It can thus be
exploited exactly as before.  Note also that~$L$ is itself the generator of
a stationary Markov process~$\hZ$, the Ornstein--Uhlenbeck diffusion.

For the purposes of normal approximation, the Stein identity 
$\ex\{\tL(f\circ g)(\tZ(0))\}=0$, for all suitable~$f$ and for~$g$ the 
component sum, as above,  yields
$$
    -\ex Lf(W) \Eq {1\over2}\ex\Blb \Bl\sn X_i^2 - 1\Br f''(W)\Brb + \ex E(f,g,X), \eqno(4.5)
$$
because of (4.4). If the right hand side of (4.5) can be unifomly bounded for
$f=f_h$ and $h\in\hh$,  a corresponding bound for the accuracy of the
normal approximation is obtained.  In particular, it is necessary to show that 
$$
   \ex \Blm \sn X_i^2 - 1\Brm \qquad{\rm and}\qquad \ex|E(f,g,X)|         \eqno(4.6)
$$
are both small.  The latter can be estimated by a Lyapounov term, if $\|f'''\|_\infty$
is finite, but requires more work if this cannot be assumed, and the former can also be 
expected to be small if the Lyapounov ratio is small.  

The concrete arguments above have always involved properties of the solution~$f_h$ to
the Stein equation $Lf = h$.  For the standard normal distribution, Stein~(1972) gave
bounds for $f_h$ and $f'_h$ in the supremum norm, when~$h$ is an indicator function.
Indeed, the explicit form of the solution~$f_h$ given in (2.3) greatly facilitates 
such calculations.  For other Stein equations, controlling the solutions can be
much more difficult.  However, for~$L$ the generator of a stationary Markov process~$\hZ$
with stationary distribution~$\law(Z)$,
the equation $Lf=h$ is the so-called Poisson equation, and the solution, under
appropriate conditions, is the recurrent potential
$$
     f_h(w) \Eq - \int_0^\infty \ex\{h(\hZ(t))\giv \hZ(0) = w\}\,dt,       \eqno(4.7)
$$
assuming as always that $\ex h(Z) = 0$.

\bsk\msk
\vfil\eject
\nin{\sl 5. Exchangeable pairs.}\msk
In his 1986 monograph `Approximate Computation of Expectations', Stein begins by introducing the
 notion of an exchangeable pair of random variables for the construction of Stein identities.
As usual, the aim is to approximate the distribution of some random element~$W$.  Suppose that~$W'$
is another random variable, defined on the same probability space as~$W$, such that $(W,W') =_d (W',W)$. 
Then $(W,W')$ is called an exchangeable pair.  Hence, if~$F$ is any antisymmetric real function
on $\ss\times\ss$, it is immediate that $\ex F(W,W') = 0$ (if the expectation exists), and hence
also that $\ex\{\ex^W F(W,W')\} = 0$, where $\ex^W$ denotes the conditional expectation given~$W$.
If~$\ss$ is finite, and connected in the sense that $w_1$ and~$w_2$ are neighbours if
$\pr[(W,W') = (w_1,w_2)] > 0$, Stein showed that every function~$g$ such that $\ex g(W) = 0$ can
be represented in the form $g(W) =  \ex^W F(W,W')$, characterizing the distribution~$\law(W)$
in terms of conditional expectations of the form $\ex^W F(W,W')$.

To convert this statement into a characterizing operator $\tL$ on~$\xx$, one needs to define
an appropriate subcollection of antisymmetric functions~$F$.  For normal approximation,
Stein proposed the functions
$$
    F(w,w') \Eq (w-w')(f(w) + f(w')),  \eqno(5.1)
$$
for $f\in\xx$.  If the linear regression condition
$$
   \ex^W W' \Eq (1-\l)W  \eqno(5.2)
$$
is satisfied for some $0 < \l < 1$ (already implying that $\ex W = 0$), 
 Stein~(1986, Theorem III.1) showed that
$$\eqalign{
   |\pr[W \le w] - \Phi(w)| &\Le 2\sqrt{\ex\Bl 1 - {1\over 2\l}\ex^W\{(W'-W)^2\}\Br^2}\cr
      &\hbox{}\qquad + {1\over (2\pi)^{1/4}} \sqrt{{1\over\l}\ex|W-W'|^3}.}   \eqno(5.3)
$$
Thus, provided that $W$ and~$W'$ are typically close to one another, and that the
conditional expectation $\ex^W\{(W'-W)^2\}$ is concentrated near the value $2\l$,
the linear regression property implies that the distribution of~$W$ is close to 
standard normal.
This theorem, and its refinements (in particular allowing the linear regression
condition to be only approximately satisfied), provide a powerful and flexible tool
for estimating the accuracy of normal approximation in a wide variety of settings,
with complicated dependence structures.

How is this theorem proved?  The basic idea is simple.  Take the antisymmetric function~$F$
from (5.1), giving
$$
   \ex \{(W-W')(f(W) + f(W'))\} \Eq 0. \eqno(5.4)
$$
This, with liberal use of (5.2), can equivalently be written as
$$\eqalign{
  \ex\{Wf(W) - f'(W)\} &\Eq \ex\Blb f'(W) \Bl {(W-W')^2 \over 2\l} - 1 \Br \Brb \cr
     &\hbox{}\qquad + {1\over 2\l}\ex\Blb (W-W')\{f(W) - f(W') - (W-W')f'(W)\}\Brb.}
   \eqno(5.5)
$$
Note that the left hand side is just $-\ex \{Lf(W)\}$, for~$L$ the Stein operator
for the standard normal distribution.
The first term on the right hand side is small if $\ex^W\{(W'-W)^2\}/2\l$ is
concentrated close to~$1$ (and (5.2) implies that $\ex\{(W'-W)^2\} = 2\l\ex W^2$,
so that having $\ex W^2$ close to~$1$ is certainly a good idea), and, for twice
differentiable functions~$f$, the second term is plausibly small
if  ${1\over\l}\ex|W-W'|^3$ is small.  The detailed argument leading to (5.3)
is somewhat more complicated, but Stein's method has already accomplished the
hard work.

The choice of antisymmetric function $F(w,w') = (w-w')(f(w) + f(w'))$ is not
the simplest: one could instead have chosen $F(w,w') = f(w') - f(w)$.  This
would in similar style give
$$\eqalign{
  0 &\Eq \ex\{f(W') - f(W)\} \cr
    &\Eq \ex\Blb(W'-W)f'(W) + {1\over2}(W'-W)^2f''(W)\Brb + \ex E'(f,W,W'),
}   \eqno(5.6)
$$
where
$$
   E'(f,W,W') \Eq f(W') - f(W) - (W'-W)f'(W) - {1\over2}(W'-W)^2 f''(W).
$$
Using the linear regression condition $\ex^W W' \Eq (1-\l)W$ of (5.2), which
as above also implies that $\ex\{(W'-W)^2\} = 2\l\ex W^2$, it follows that
$$
   0 \Eq -\l \ex\{-Wf'(W) + f''(W)\}
     + \l\ex \Blb \Bl {\ex^W(W'-W)^2 \over 2\l} - 1\Br f''(W)\Brb
     + \ex E'(f,W,W').  \eqno(5.7)
$$
Cancelling the factor~$\l$, this gives a close parallel to (5.5), but with
$f'$ replacing~$f$ in the analogy, as was the case when using the
generator approach in Section~4.  Note that, in this argument, exchangeability
is not used: $\ex\{f(W')-f(W)\}=0$ requires only that $W=_d W'$ (R\"ollin, 2008).

In many examples, the exchangeable
pair can be realized as a pair of successive states in a stationary
{\it reversible\/} Markov chain.  In these circumstances, (5.6) can also be interpreted
as the Stein identity from the generator approach.  This is indeed the
case for the sums of independent random variables considered in the
previous section.  The linear regression condition is satisfied with
$\l = 1/n$, and expression (4.5) is an exact parallel to (5.7).

The approach using exchangeable pairs has proved extremely useful in many contexts.
In the book Diaconis and Holmes~(2004), for instance, there are chapters showing 
how Stein's method and exchangeable pairs can be effectively exploited in a 
variety of quite disparate 
settings.  They are used, amongst other things, to reduce variance in simulation experiments
(Stein, Diaconis, Holmes \& Reinert 2004), to
estimate rates of mixing for some ergodic Markov chains (Diaconis 2004), to prove a central limit
theorem for the number of descents in a random permutation (Fulman 2004), and to establish the
large sample properties of the bootstrap, including a variant of the
blockwise bootstrap for some dependent samples (Holmes \& Reinert 2004). 
Other examples that have since been influential include those in Rinott \& Rotar~(1997).
The combination of linear regression condition and the exchangeable pair has
recently been effectively developed in the setting of multivariate normal
approximation by Reinert \& R\"ollin~(2009).

\bsk\msk

\nin{\sl 6. The concentration inequality approach.}\msk
The Berry--Esseen theorem is framed in terms of the Kolmogorov distance, in which
the class~$\hh$ of test functions~$h$ consists of the indicators of half
lines.  The corresponding functions~$f_h$ defined by (2.3) do not have 
$\|f_h''\|_\infty < \infty$, since their derivatives have a jump.
This means that the simple argument given in Section~2 (c.f.\ Erickson 1974) 
cannot be used directly
to prove normal approximation with respect to this distance.
However, the effect of the jump on the expectations to be bounded can be
controlled, if it is known that the concentration function of the distribution
of~$W$ is suitably bounded.  This was an essential step in Stein's original proof 
of the Berry--Esseen 
theorem for sums of independent and identically distributed random variables. 
Astonishingly, his proof of the concentration inequality also
made use of his method in a most ingenious way. 
We now reproduce his proof in detail, in the very simplest context.  

Let
$X_1 , X_2 , \ldots , X_n$ be a sequence of independent and identically distributed 
random variables such that $\ex X_1 = 0$, $\ex X_1^2 = 1/n$ and $\ex |X_1|^3 < \infty$; 
set $\b = n^{3/2}\ex |X_1|^3$. The Berry--Esseen theorem states that, 
for every real~$z$,
$$
                  |F_n(z) - \Phi(z)| \Le C\b n^{-1/2},            \eqno(6.1)
$$
where $F_n$ is the distribution function of $\sn X_i$ , $\Phi$ is the standard normal 
distribution function, and $C$ is an absolute constant.
                                                                       
    Stein's proof, to be found in Ho \& Chen~(1978), proceeds as follows. 
Let $W_n =\sn X_i$ , $W_{n-1} = \sum_{i=1}^{n-1} X_i$. 
Then, for any absolutely continuous function~$f$ with bounded derivative, we have
$$
       \ex \{W_n f (W_n )\} \Eq \sn \ex \Blb X_i f\Bl \sum_{j\ne i} X_j + X_i \Br\Brb, \eqno(6.2)
$$                                        
which by independence, symmetry and the fact that $\ex X_n = 0$, gives
$$\eqalign{
    \ex \{W_n f (W_n )\} &\Eq   n\ex [X_n f (W_{n-1} + X_n )] 
         \Eq n\ex [X_n (f (W_{n-1} + X_n ) - f (W_{n-1}))] \cr
     &\Eq n\ex \Blb X_n \int_0^{X_n} f'(W_{n-1} + t)\,dt \Brb \cr
     &\Eq \ex \int_{-\infty}^\infty  f'(W_{n-1} + t)K(t)\,dt ,} \eqno(6.3)
$$      
where
$$
           K(t) \Eq n\ex X_n [I(0 < t < X_n ) - I(X_n < t < 0)] .\eqno(6.4)
$$
Hence we obtain the identity
$$                                              
     \ex \{W_n f (W_n )\} \Eq \ex \int_{-\infty}^\infty f'(W_{n-1} + t)K(t)\,dt . \eqno(6.5)
$$                                           
Now the function $K(t)$ is nonnegative, satisfying $K(-\infty) = 0$ and $K(\infty) =
0$, and, by setting $f (w) = w$, (6.5) yields
$$                                                   
                     \int K(t)\,dt \Eq \ex W_n^2 \Eq 1 ,             \eqno(6.6)
$$
showing that $K$ is in fact a probability density function; we note also, for future
reference, that  
$$
  \int |t| K(t)\,dt \Eq \b/2\sqrt n.  \eqno(6.7)
$$
Hence, from (6.5), the Stein operator for standard normal approximation applied to~$W_n$
can be expressed as 
$$                                              
   \ex \{f'(W_n) - W_n f (W_n )\} 
    \Eq \ex \int_{-\infty}^\infty \{f'(W_{n-1} + X_n) - f'(W_{n-1} + t)\}K(t)\,dt , \eqno(6.8)
$$ 
and it remains only to show that the expression on the right hand side is small.  

Clearly,
as before, there is a very simple argument to do so if~$f_h''$ is uniformly bounded for
our test functions~$h$, but $f'_h$ has a jump of size~$1$ at~$a$ when $h = \bone_{(-\infty,a]}$,
and its derivative may be big close to~$a$ if~$a$ is large.  Thus it is necessary to show 
that the random variable~$W_{n-1}$ rarely takes a value close enough to~$a$ to make the
expected difference $\ex^{W_{n-1}}|f'(W_{n-1} + X_n) - f'(W_{n-1} + t)|$ appreciable.
Note that, from (6.7), the relevant values of~$t$ here are of order $1/\sqn$, and this
is also the scale appropriate to values of~$X_n$; hence it is necessary to be able to
show that~$W_{n-1}$ has suitably small probability of being close to~$a$ on this scale.
It turns out that the right hand side of (6.8) can indeed be shown to be of 
the required order $O(\b/\sqn)$, by using the following concentration inequality.

\proclaim Lemma. For all real $a$ and $b$ such that $a < b$, 
$$
                        \pr(a \le W_{n-1} \le b) \Le b - a + 2\b/\sqrt n .
$$

\proof
We begin with a simple inequality. By (6.7), we have
$$                                                                             
   \int_{t>\b/\sqn}K(t)\,dt \Le (n^{1/2}/\b)  \int_{t>\b/\sqn}|t|K(t)\,dt 
              \Le (n^{1/2}/\b) \int  |t|K(t)\,dt \Eq {1\over2}.
$$
    This and (6.6) yield
$$
     \int_{t\le\b/\sqn}K(t)\,dt  \Eq \int K(t)\,dt - \int_{t>\b/\sqn}K(t)\,dt 
      \ \ge\ {1\over2} .  \eqno(6.9)
$$ 

Now let $a$ and $b$ be two real numbers such that $a < b$. For any real
$x > 0$, define
$$                         
   g_x (w) \Eq \cases{ -{1\over2} (b - a) - x &if $w \le a - x$\cr
                         w - {1\over2}(a + b)  &if $a - x \le w \le b + x$\cr
                        {1\over2}(b - a) + x  &if $b + x \le w$.} \eqno(6.10)
$$ 
Clearly, $g_x$ is an absolutely continuous function, and $g_x'(w) = I(a-x \le w \le b+x)$.
Taking $f=g_{\b/\sqn}$, we have
$$\eqalign{                                       
  \ex  \int f'(W_{n-1} + t)K(t)\,dt &\Eq \ex \int I(a - \b/\sqn \le W_{n-1}+t \le b + \b/\sqn)K(t)\,dt \cr
           &\ \ge\ \ex  I(a \le W_{n-1} \le b)\int I(|t| \le \b/\sqn)K(t)\,dt.} \eqno(6.11)
$$
By (6.9), the right hand side is at least $ {1\over2}\pr(a \le W_{n-1} \le b)$.
But now, in view of (6.5), this gives
$$
      \pr(a \le W_{n-1} \le b) \Le 2\ex W_n f (W_n ) \Le 2\ex |W_n f (W_n )| \Le b - a + 2\b/\sqn,
      \eqno(6.12)
$$                         
because $\|f\|_\infty \le {1\over2} (b - a) + \b/\sqn$ from (6.10),
and $\ex |W_n | \le (\ex W_n^2)^{1/2} = 1$.
\ep

    Stein's idea of concentration inequality was extended beyond the case of i.i.d.\ 
random variables in
Chen~(1986). For locally dependent and not necessarily identically distributed random
variables $X_1 , X_2 , \ldots , X_n$ and a random vector~$\z$ depending on a relatively small
subset of $\{X_1 , X_2 , \ldots , X_n\}$, a randomized concentration inequality can be proved
for $P^{\z}(a\z \le W \le b\z)$.  This is used to derive good bounds on the
difference with respect to Kolmogorov distance between the distribution
the standardized sum~$W$ of the $X_i$'s
and the standard normal distribution, thereby extending the classical Berry--Esseen theorem. 
Non-uniform concentration inequalities were also developed in Chen \& Shao (2001, 2004), 
and were used to obtain a non-uniform Berry--Esseen theorem.

    If a random variable $T$ can be expressed as $T = W + \D$, where the Kolmogorov
distance between $\law(W)$ and $\law(Z)$ is known to be small and where~$\D$ is a 
relatively small
random variable which may depend on~$W$, a randomized concentration inequality
for $\pr(z \le W \le z + |\D|)$ again gives good bound for the Kolmogorov distance
between $\law(W)$ and $\law(Z)$. This idea
was used in Chen \& Shao~(2007) to prove normal approximation for nonlinear
statistics, and in Barbour \& Chen~(2005) for the permutation distribution of matrix
correlation statistics.
    
\bsk\msk
\nin{\sl 7. Poisson approximation.}\msk
One of the biggest successes of Stein's method has been in Poisson approximation.
The classical limit theorem, that $\Bi(n,\l/n) \to \Po(\l)$ is taught in most
first courses on probability theory, and the proof uses explicit computation
of the point probabilities.  It was only rather recently (Prohorov, 1953) that
a good bound on the difference between the two distributions was given,
in the form
$$
    \dtv(\Bi(n,p),\Po(np)) \Le cp\min\{1,np\},  \eqno(7.1)
$$
for an explicit constant~$c$, where the total variation distribution~$\dtv(P,Q)$
between two probability measures is the supremum of $|P(A)-Q(A)|$ over measurable
sets~$A$.  Note that, if $W$ and~$Z$ are random variables on some~$\ss$, then
$$
   \dtv(\law(W),\law(Z)) \Eq {1\over2} \sup_{h\in\hh_0} |\ex h(W) - \ex h(Z)|,  \eqno(7.2)
$$
where $\hh_0$ is the set of measurable real functions~$h$ on~$\ss$ with
$\|h\|_\infty \le 1$.  This representation suggests that Stein's method may
be well suited to total variation approximation: its development in the case
of Poisson approximation is due to Chen~(1975).

The first step is to derive a suitable Stein operator on real valued functions
on the non-negative integers~$\integ_+$.  By analogy with the derivation of the
Stein operator for standard normal approximation in Section~2, letting~$p(j)$
denote the Poisson~$\Po(\l)$ probability of the point~$j$, we have
$$\eqalign{
  0 &\Eq \sum_{j=0}^\infty \D_*\{f(j)p(j)\}  \cr
    &\Eq \sum_{j=0}^\infty  \{f(j)  - jf(j-1)/\l\}p(j)
}  \eqno(7.3)$$
for all $f\in\ff$, the set of real functions~$f$ on~$\integ_+$ such that 
$f(j)p(j) \to 0$ as $j\to\infty$,
where $\D_* g(j) = g(j)-g(j-1)$ denotes the backward difference operator,
and $\D_*\{f(0)g(0)\} = f(0)p(0)$.
Thus any function~$h$ of the form $h(j) =  f(j) - jf(j-1)/\l$ with $f\in\ff$
automatically has $\ex h(Z)=0$, if $Z \sim \Po(\l)$.  Conversely, for any~$h$
such that $\ex |h(Z)| < \infty$, the function $h(j) - \ex h(Z)$ can be written
in the form $ f(j) - jf(j-1)/\l$, with $f\in\ff$.  This is because, 
if $g(j) =  f(j) - jf(j-1)/\l$, then
$$\eqalign{
  \sum_{j=0}^J g(j)p(j) 
   &\Eq \sum_{j=0}^J \D_*\{f(j)p(j)\} \Eq f(J)p(J), 
} \eqno(7.4)$$
and, in reverse, defining $f$ by $f(J)p(J) = \sum_{j=0}^J g(j)p(j)$ for any~$g$
such that $\ex |g(Z)| < \infty$ yields $ f(j) - jf(j-1)/\l = g(j)$.  
Hence we can take~$f=f_h$ to be defined by
$$\eqalign{
   f_h(J)p(J) 
            &\Eq \sum_{j=0}^J  \{h(j)  - \ex h(Z)\}p(j) \cr
            &\Eq - \sum_{j=J+1}^{\infty} \{h(j) - \ex h(Z)\}p(j)}, \eqno(7.5)
$$
noting that it is then directly checked that $f_h\in\ff$.  
As before, this allows us to write
$$
   \ex h(W) - \ex h(Z) \Eq \ex\{f_h(W) - W f_h(W-1)\},       \eqno(7.6)
$$
for any random variable~$W$ on~$\integ_+$ for which the expectations exist;
the value of $f_h(-1)$ can be chosen arbitrarily, since it plays no role.
Thus we are led to a Stein operator for the Poisson distribution, usually
expressed, by shifting the argument of~$f$ and multiplying through by~$\l$,
in the equivalent form
$$
    Lf(w) \Eq \l f(w+1) - w f(w),                               \eqno(7.7)
$$
for functions $f$ such that $f(j+1)p(j) \to 0$ as $j\to\infty$.

To illustrate its use, let $X_1,\ldots,X_n$ be independent Bernoulli random
variables, with $X_i \sim \Be(p_i)$, and let $W = \sn X_i$.  If the~$p_i$
are not too large, Le Cam~(1960) gave an analogue of Prohorov's (1953)
approximation for the total variation distance between $\law(W)$ and~$\law(Z)$, 
where $Z \sim \Po(\l)$ and $\l = \sn p_i$: he showed that
$$
    \dtv(\law(W),\law(Z)) \Le 8\min(1,\l^{-1})\sn p_i^2,          \eqno(7.8)
$$
if $\max_i p_i \le 1/4$, using an operator approach that relies heavily on
independence.  To apply Stein's method, we need to bound the right hand side
of~(7.6).  As for normal approximation, the term $\ex\{Wf_h(W)\}$ breaks up
into elements of the form $\ex\{X_if_h(W)\}$, and since $X_if(W) = X_if(W_i+1)$
a.s., where $W_i = W-X_i$, it follows by the independence of $X_i$ and~$W_i$
that
$$
    \ex\{X_i f_h(W)\} \Eq p_i\ex f_h(W_i+1);\qquad \ex\{Wf_h(W)\} 
      \Eq \sn p_i\ex f_h(W_i+1).
$$
This is to be compared with 
$$
    \l \ex f_h(W+1) \Eq \sn p_i \ex f_h(W_i + X_i + 1).
$$
Taking the difference, it follows easily that
$$
  |\ex h(W) - \ex h(Z)| \Le \sn p_i\ex |f_h(W_i+X_i+1) - f_h(W_i+1)|
                        \Le \sn p_i^2 \|\D f_h\|_\infty,           \eqno(7.9)
$$
where $\D f(j) = f(j+1)-f(j)$ denotes the forward difference.
Since $\sup_{h\in\hh_0}\|\D f_h\|_\infty \le 2\min(1,\l^{-1})$, it follows 
immediately that
$$
    \dtv(\law(W),\law(Z)) \Le \min(1,\l^{-1})\sn p_i^2,                  \eqno(7.10)
$$
improving the constant in Le Cam's bound, and removing the restriction on the
maximum value of the~$p_i$'s.  The argument was also refined to establish
a {\it lower\/} bound for the error, having exactly the same order (Barbour \&
Hall, 1984).

Of course, the real advance of the method in this context is the ease with which
dependence can be accommodated.  For instance, if $W$ can be written as $X_i + Y_i + \tW_i$,
with $Y_i$ a random variable in~$\integ_+$ and $\tW_i$ almost independent of $X_i$, 
then one can write
$$\eqalign{
    \ex\{X_i f_h(W)\} &\Eq \ex\{X_i f_h(1+Y_i+\tW_i)\} \cr
      &\Eq \ex\{X_i f_h(1+\tW_i)\}   +  \ex\{X_i (f_h(1+Y_i+\tW_i)  -  f_h(1+\tW_i))\},
                                                                           }\eqno(7.11)
$$
with the first term close to $p_i \ex f_h(\tW_i + 1)$ and the second bounded by
$\ex (X_iY_i)\|\D f_h\|_\infty$.  Comparing this to
$$
   p_i \ex f_h(W+1) \Eq p_i \ex f_h(X_i + Y_i + \tW_i),
$$
it follows immediately that
$$\eqalign{
    &|\ex h(W) - \ex h(Z)| \cr
   &\quad\Le \sn\{p_i^2 + p_i\ex Y_i + \ex(X_iY_i)\}\min\{1,\l^{-1}\}
    + \sn \ex|\ex(X_i\giv \tW_i) - p_i|\,\|f_h\|_\infty,   }                \eqno(7.12)
$$
and $\sup_{h\in\hh_0}\|f_h\|_\infty \le 2\min(1,\l^{-1/2})$.  This leads to useful
Poisson approximation in total variation for a wide variety of sums of dependent 
Bernoulli random variables.  

Poisson approximation for the sum of independent Bernoulli random variables 
could also have been approached by
deriving a Stein identity for~$W$, and then deducing the Stein operator appropiate
to the problem.  A stationary Markov chain~$\tZ$ with $\tZ(0) = (X_1,\ldots,X_n)$ 
can be constructed exactly as in Section~4, and its generator is
$$
   \tL f(z) \Eq {1\over n}\sn \ex\{f(z + \be_i(X_i-z_i)) - f(z)\}; \eqno(7.13)
$$
specializing to functions $f\circ g$ with $g(z) = \sn z_i$, this gives,
after some manipulation,
$$\eqalign{
   &n\tL (f\circ g)(z) \cr
   &\Eq \l\{f(g(z)+1) - f(g(z))\} - g(z)\{f(g(z)) - f(g(z)-1)\} 
     +  \sn z_i p_i (\D^2f)(g(z)-1).}                          \eqno(7.14)
$$
Thus, writing $g(z)=w$, and dropping the factor~$n$, this expression is close to $Lf(w)$, 
where
$$
    Lf(w) \Eq \l \D f(w) - w\D f(w-1) ,                                 \eqno(7.15)
$$
the same generator as in (7.7), except that now $\D f(w-1)$ plays the part of $f(w)$.
Note that~$L$ is the generator of a well-known stationary, reversible Markov process,
the immigration--death process with immigration rate~$\l$ and with unit {\it per capita\/}
death rate.  Poisson approximation for sums of independent Bernoulli random variables
can now be deduced directly from (7.14), since it gives
$$
   |\ex Lf_h(W)| \Le \sn p_i^2 \|\D^2 f_h\|_\infty, \eqno(7.16)
$$
which is precisely the bound given in (7.9).  

In more complicated circumstances, with
dependent summands, it may still sometimes be easier to 
start with the version of~$Lf$ given in (7.7) or in (7.15) and to make direct calculations,
rather than to set up a Stein identity for~$W$ as the first step.
On the other hand, one step of the chain~$\tZ$ in the independent setting would 
yield an exchangeable pair, and constructing an exchangeable pair also offers many possibilities
in dependent settings in which direct calculation may seem unattractive. 

The identification of the distribution $\Po(\l)$ as that of a stationary immigration--death
chain can be very simply extended to Poisson processes. A Poisson point process~$\Pi$
with intensity measure~$\l$
on a space~$\uu$ has distribution~$\pp(\l)$, the equilibrium distribution of a 
spatial immigration--death process~$\tPi$
with immigration intensity measure~$\l$ over~$\uu$ and with unit {\it per capita\/} death
rate.  Letting~$\d_u$ denote the point mass at~$u$, $\xi$ a Radon measure on~$\uu$
and~$f$ a bounded real function on the space of such Radon measures, the generator
of the process~$\tPi$ can be expressed as
$$
   Lf(\xi) \Eq \int_{\uu}\{f(\xi + \d_u) - f(\xi)\}\l(du) + 
                   \int_{\uu} \{f(\xi - \d_u) - f(\xi)\}\xi(du).                  \eqno(7.17)
$$
The corresponding Stein equation $Lf = h$, for bounded functions~$h$ for
which $\ex h(\Pi) = 0$, can be solved in the form
$$
   f_h(\xi) \Eq -\int_0^\infty \ex\{h(\tPi(t)) \giv \tPi(0) = \xi\}\,dt,             \eqno(7.18)
$$
as an example of (A4.9), and coupling arguments can be used to bound the first
and second differences of~$f_h$.  However, at least for the class of test functions appropriate
to total variation approximation, the bounds on these differences do not in general decrease
with increasing $\l(\uu)$, in contrast to the case of the Poisson distribution~$\Po(\l)$,
which corresponds to~$\uu$ consisting of a single point.  Much more on this topic
can be found in Barbour \& Brown~(1992) and in \BHJ~(1992, Chapter~X).

The simplest application is once again in the independent setting.  With $X_1,\ldots,X_n$
as defined earlier in the section, let $\uu=\{1,2,\ldots,n\}$, and set $\Xi = \sn X_i\d_i$.
To see how close the distribution of~$\Xi$ is to that of a Poisson process~$\Pi$ with intensity
$\l = \sn p_i\d_i$, we need to examine $\ex Lf(\Xi)$.
Now, much as for the random variable case,
$$\eqalign{
   \ex\Blb \int_{\uu} \{f(\Xi - \d_u) - f(\Xi)\}\Xi(du) \Brb
    &\Eq -\sn \ex (X_i\{f(\Xi_i + \d_i) - f(\Xi_i)\})  }
$$
where $\Xi_i = \Xi - X_i\d_i$ is independent of~$X_i$, and
$$
   \ex\Blb \int_{\uu} \{f(\Xi + \d_u) - f(\Xi)\}\l(du) \Brb
     \Eq \sn p_i \ex \{f(\Xi + \d_i) - f(\Xi)\}.
$$
Adding the two yields
$$
   \ex Lf(\Xi) \Eq \sn p_i\{f(\Xi_i + \d_i + X_i\d_i) - f(\Xi_i+X_i\d_i) 
            - f(\Xi_i + \d_i) + f(\Xi_i)\}, \eqno(7.19)
$$
giving an immediate bound
$$
   |\ex h(\Xi) - \ex h(\Pi)| \Eq
   |\ex Lf_h(\Xi)| \Le \sn p_i^2 \sup_{u,\xi}|f_h(\xi + 2\d_u) - 2f_h(\xi+\d_u) + f_h(\xi)|,
                                                                                   \eqno(7.20)
$$
and the second difference can be bounded starting from (7.18); for total variation,
the supremum is just~$2$.  Note that the argument is in every sense equivalent to that
for random variables.

For dependent random variables~$X_i$, there is again a parallel.  Suppose now that, for each~$i$,
we can write $\Xi = X_i\d_i + H_i + \tXi_i$, where $H_i$ is a Radon measure and $X_i$
and $\tXi_i$ are almost independent.  Then
$$\eqalign{
  &\ex\Blb \int_{\uu} \{f(\Xi - \d_u) - f(\Xi)\}\Xi(du) \Brb 
    \Eq -\sn \ex (X_i \{f(\tXi_i + \d_i)  - f(\tXi_i)\}) \cr
      &\qquad\qquad\qquad\hbox{}   - \sn\ex(X_i\{f(H_i + \tXi_i) - f(H_i + \tXi + \d_i)  
             - f(\tXi_i + \d_i)  + f(\tXi_i)\}), 
} \eqno(7.21)$$
with the first term, much as for (7.11), close to 
$\sn p_i \ex \D_i f_h(\tXi_i)$, and the second bounded by
$\ex (X_i Y_i)\sup_{i,j}\|\D_{ij} f_h\|_\infty$: here, $\D_i$ denotes the first forward
difference in the direction~$i$, and $\D_{ij}$ the corresponding second differences, and
$Y_i = H_i(\uu)$. This is to be added to
$$\eqalign{
  & \ex\Blb \int_{\uu} \{f(\Xi + \d_u) - f(\Xi)\}\l(du) \Brb  \cr
     &\qquad\Eq \sn p_i \ex \{f(X_i\d_i + H_i + \tXi_i + \d_i) - f(X_i\d_i + H_i + \tXi_i)\}\cr
     &\qquad\Eq \sn p_i\ex  \D_i f(\tXi_i) 
       + \sn p_i \ex \{\D_i f(X_i\d_i + H_i + \tXi_i) - \D_i f(\tXi_i)\},
} \eqno(7.22)$$
whose second term is bounded by
$$
    \sn p_i\{p_i + \ex Y_i\}\sup_{i,j}\|\D_{ij} f_h\|_\infty.  \eqno(7.23)
$$
Thus, in parallel to (7.12), we obtain
$$
 |\ex h(\Xi) - \ex h(\Pi)| \Le \sn\{p_i^2 + p_i\ex Y_i + \ex(X_iY_i)\}c_1(h)
    + \sn \ex|\ex(X_i\giv \tXi_i) - p_i|\,c_0(h);                              \eqno(7.24)
$$ 
again, in contrast with the situation for random variables, the best general
bounds for $c_0(h)$ and~$c_1(h)$ for functions~$h$ with $\|h\|_\infty \le 1$
are each~$2$, so that there is in general no decrease when $\l(\uu)$ increases. 

The total variation bound for the approximation of~$\law(\Xi)$ by $\PP(\l)$ implied 
by~(7.24) was demonstrated in a wide variety of examples in \AGG~(1989, 1990).
They also used the bound as a basis for compound Poisson approximation, by showing
that random variables that might be expected to have approximately such a distribution
can frequently be represented as functions of Poisson processes --- a clump of
observations of size~$j$ at point~$x$ is associated with the point mass~$\d_{(x,j)}$,
and~$\uu$ consists of the set of such pairs.  Their technique turns out to be
widely applicable in practice.  Compound Poisson approximation can also be
approached directly, using a Stein operator derived in \BCL~(1992) and the
method for exploiting it derived in Roos~(1994).  For this
operator, there are difficulties in bounding the solutions of the corresponding
Stein equation --- see Barbour \& Utev (1998, 1999) and Barbour \& Xia~(1999).
The expository paper of Barbour \& Chryssaphinou~(2001) gives
a discussion of the various approaches to compound Poisson approximation using 
Stein's method, and a more detailed overview of Poisson approximation can also
be found in Barbour~(2001).

\bsk
\nin{\sl 8. Ramifications.}\msk  
This article is intended more as a personal view of the extraordinary impact that
Stein's method has had, rather than as an exhaustive survey; the subject is just too
large.  Much more information is to be
found in the books and expository articles that have already been written on
Stein's method.  The primary source for the general method is in those by Stein 
himself (1986, 1992).  For Poisson approximation and related themes, the papers
of \AGG~(1989, 1990) and the book by \BHJ~(1992) are good starting points. 
The Diaconis \& Holmes~(2004) collection has lots of ideas concerned with the
use of exchangeable pairs.  The tutorial lectures in Barbour \& Chen (2005a)
also furnish a digestible introduction to many aspects of Stein's method,
and the accompanying workshop proceedings (Barbour \& Chen 2005b) illustrate how wide the scope
of Stein's method has become.

This latter aspect, the enormous breadth of research that has developed out of
Stein's ideas, should nonetheless be briefly emphasized.  For instance, the
method has been applied in the context of approximation by distributions unrelated
to either Poisson or normal.  Luk~(1994) studied approximation by Gamma distributions,
Loh~(1992) by the multinomial distribution, Pek\"oz~(1996) by the geometric distribution,
Brown \& Xia~(2001) by the equilibrium distribution of a birth and death process,
and G\"otze \& Tikhomirov~(2005) by the Wigner semi-circle law.  This latter paper
is but one of a number of applications of Stein's method in random matrix theory:
for example, there is the recent work of Meckes (2008) giving sharp rates for 
normal approximation to the distribution of linear functions of random orthogonal 
and unitary matrices,  Fulman's~(2009)
normal approximation to the distribution of the trace of a random matrix from a 
compact Lie group, and Chatterjee's~(2009) results discussed below.  
Gaussian and Poisson process approximation by Stein's method have 
already been briefly mentioned; Reinert~(1995) showed how to use the method to
prove weak laws of large numbers for empirical processes.

Another fascinating development coming from the original normal and Poisson
applications is the relationship between the method and certain biasing
constructions.  The `coupling' approach to Poisson approximation can be 
interpreted in terms of the
size biasing characterization of the Poisson distribution.  Size biasing
can also be naturally introduced into Stein's method for the normal distribution,
when the random variables under consideration are non-negative (Goldstein \& Rinott, 1996).
In Goldstein \& Reinert~(1997, 2005), a new biasing construction, known as zero biasing,
was introduced. The construction associates a zero biased distribution to any
distribution with zero mean and finite variance, in a way that fits very neatly
with the Stein operators for normal approximation in one and higher dimensions.

An area in which Stein's method has proved extremely fruitful
is that of random geometrical graphs, in which points in space are a
taken to be neighbours if they are close in an appropriate sense.  Such models
arise very naturally in many branches of statistics.  Starting from a paper
of Avram \& Bertsimas~(1993), itself based on Baldi \& Rinott~(1989), Stein's
method has been used in a systematic way to prove many
theorems about normal approximation to the characteristics of such graphs; see,
for example, Penrose \& Yukich~(2005) for more details. 

There has recently been an explosion of interest in the application of Stein's
method to functionals of Gaussian processes and fields.  Using the method,
combined with Malliavin calculus, it has been possible to generalize, refine
and unify many central and non-central limit theorems for multiple Wiener--It\^o integrals.
A first step in this direction is to be found in Nourdin \& Peccati~(2009a),
sharpened and refined in Nourdin \& Peccati~(2009b) to obtain one term Edgeworth
expansions.  Their approach has a wide range of application, and it is too
soon for a complete picture of the possible developments to be discerned.

Stein's (1986) monograph contains a chapter on large deviations, which can be
seen as work in progress.  This project has now been successfully developed
by Chatterjee~(2007a), who is able to prove strong concentration of measure inequalities
(to distinguish them from bounds on the concentration function) using the Stein
approach, in a variety of interesting settings.  This is but one of his recent 
 impressive results.  In Chatterjee~(2009), he presents a rather general 
technique for proving central limit approximations for the distributions of 
linear statistics derived from high dimensional random matrices,
using Stein's method and a notion of second order Poincar\'e inequalities.
One of the classes of matrices that he studies is that of random Gaussian Toeplitz matrices; for these 
matrices, a central limit theorem involving the spectrum is proved, even 
though the limiting formula for the associated variance is not known. 
For real-valued functionals $f(X)$ of independent random 
variables $X=(X_1,X_2,\dots,X_n)$, an abstract theorem bounding 
the distance between the distribution of $f(X)$ and the appropriate normal distribution 
is proved in Chatterjee~(2008). The main part of the bound,  a discrete version of his
second-order Poincar\'e inequality, is expressed in terms 
of the variance of a random variable~$T$, which is constructed using an independent copy $X'$ 
of~$X$. In the case of `local' functions~$f$, this variance is shown to be relatively
accessible, and the approximation to have a number of interesting applications.
He has also proved an approximation, with error bounds, to the distribution of the local field
in the high temperature phase of the Sherrington--Kirkpatrick spin glass model: the
approximation is by means of a two part mixture of normal distributions (Chatterjee, 2010).

It is not only the method itself that is broad in its scope:
the range of application of Stein's method is equally impressive.  Stein~(1986)
gives applications to counting Latin rectangles, to random allocations, to the
binary expansion of a random integer and to isolated trees in a Bernoulli random
graph.  Other applications include the analysis of molecular sequences
(Arratia, Gordon \& Waterman 1990; Neuhauser 1994), extreme value theory (Smith 1988),
reliability theory (Godbole 1993), random fields (Takahata 1983), card shuffling
(Fulman 2005), the eigenfunctions of the Laplacian on a manifold (Meckes 2009),
logarithmic combinatorial structures (Arratia, Barbour \& Tavar\'e 2003) 
and scan statistics (Glaz, Naus, Roos \& Wallenstein 1994), among many others.

Stein's method has emerged as a flexible and powerful tool for proving probability
approximations, and has stimulated research in many different directions.  
Diaconis \& Holmes observe that `for all these virtues, it still seems impossible
to give a brief, understandable explanation of the essence of Stein's method'.
We entirely agree.

\bsk
\vfil\eject
\def\refs{\hangindent10pt\hangafter=1}
\nin{\sl References.}

\smallskip
\nin
\refs
R.\ Arratia, A.\ D.\ Barbour \& S.\ Tavar\'e (2003)
  {\it Logaritnmic combinatorial structures.\/}\break
  European Math.\ Soc.\ Press, Z\"urich.

\smallskip
\nin\refs
R.\ Arratia, L.\ Goldstein \& L.\ Gordon (1989) 
  Two moments suffice for Poisson approximations: the Chen--Stein method.
  {\it Ann.\ Probab.\/}, {\bf 17}, 9--25.

\smallskip
\nin \refs
R.\ Arratia, L.\ Goldstein \& L.\ Gordon (1990) 
    Poisson approximation and the Chen--Stein method. {\it Stat.\ Science\/}
    {\bf 5}, 403--434.   

\smallskip
\nin\refs
R.\ Arratia, L.\ Gordon \& M.\ S.\ Waterman (1990) 
  The Erd\H os--R\'enyi law in distribution, for coin tossing and sequence
      matching. {\it Annals of Statistics}, {\bf 18}, 539--70.

\smallskip
\nin\refs
F.\ Avram \& D.\ Bertsimas (1993) 
On central limit theorems in geometrical probability.
{\it Ann.\ Appl.\ Probab.\/} {\bf 3}, 1033--1046.

\smallskip
\nin\refs
P.\ Baldi \& Y.\ Rinott (1989)
On normal approximations of distributions in terms of dependency graphs.
{\it Ann.\ Probab.\/} {\bf 17}, 1646--1650.

\smallskip
\nin\refs
P.\ Baldi, Y.\ Rinott \& C.\ Stein (1989) 
  A normal approximation for the number of local maxima of a random function on a graph.
     In: {\it Probability, Statistics
     and Mathematics,\/} Eds T.\ W.\ Anderson, K.\ B.\ Athreya and D.\ L.\ 
     Iglehart. Academic Press, New York, pp.~59--81.

\smallskip
\nin\refs
A.\ D.\ Barbour (1988)
  Stein's method and Poisson process convergence. 
  {\it J.\ Appl.\ Probab.\/} {\bf 25(A)}, 175--184.

\smallskip
\nin\refs
A.\ D.\ Barbour (1990)
  Stein's method for diffusion approximations.  {\it Prob.\ Theory Rel.\ Fields\/}
  {\bf 84}, 297--322.

\smallskip
\nin\refs
A.\ D.\ Barbour (2001)  Topics in Poisson approximation. In: {\it Handbook of Statistics\/}
  {\bf 19}, Eds D.\ N.\ Shanbhag and C.\ R.\ Rao, pp.~79--115.

\smallskip
\nin \refs
A.\ D.\ Barbour and T.\ C.\ Brown (1992) 
  Stein's method and point process approximation.
    {\it Stoch.\ Procs Applics\/} {\bf 43}, 9--31.

\smallskip
\nin\refs
A.\ D.\ Barbour \& L.\ H.\ Y.\ Chen (Eds.)  (2005a)
{\it An introduction to Stein's method.\/}
IMS Lecture Note Series Volume 4, World Scientific Press, Singapore.

\smallskip
\nin\refs
A.\ D.\ Barbour \& L.\ H.\ Y.\ Chen (Eds.)  (2005b)
{\it Stein's method and applications.\/}
IMS Lecture Note Series Volume 5, World Scientific Press, Singapore.

\smallskip
\nin\refs
A.\ D.\ Barbour \& L.\ H.\ Y.\ Chen (2005c)
The permutation distribution of matrix correlation statistics.
In: {\it Stein's method and applications,\/} Eds.\ A.\ D.\ Barbour \& L.\ H.\ Y.\ Chen,
IMS Lecture Note Series Volume 5, World Scientific Press, Singapore, pp.~223--246.

\smallskip
\nin\refs
 A.\ D.\ Barbour, L.\ H.\ Y.\ Chen and W.--L.\ Loh (1992) 
    Compound Poisson approximation for nonnegative random variables via Stein's method. 
    {\it Ann.\ Probab.\/} {\bf 20}, 1843--1866.

\smallskip
\nin\refs
A.\ D.\ Barbour \& O.\ Chryssaphinou  (2001)
  Compound Poisson approximation: a user's guide.  
  {\it Ann.\ Appl.\ Probab.\/} {\bf 11}, 964--1002.

\smallskip
\nin\refs
A.\ D.\ Barbour \& P.\ Hall(1984)
 On the rate of Poisson convergence. 
  {\it Math.\ Proc.\ Cam.\ Phil.\ Soc.\/} {\bf 95},  473--480.

\smallskip
\nin\refs
A.\ D.\ Barbour, L.\ Holst and S.\ Janson (1992)
    {\it Poisson Approximation.\/}  Oxford University Press.
     
\smallskip
\nin\refs
A.\ D.\ Barbour \& S.\ Utev (1998) 
  Solving the Stein Equation in compound Poisson approximation.
      {\it Adv.\ Appl.\ Probab.\/} {\bf 30}, 449--475.

\smallskip
\nin\refs
A.\ D.\ Barbour \& S.\ Utev (1999) 
      Compound Poisson approximation in total variation.
      {\it Stoch.\ Procs.\ Applics.\/} {\bf 82}, 89--125.

\smallskip
\nin\refs
A.\ D.\ Barbour \& A.\ Xia (1999)
 Poisson perturbations.  {\it ESAIM: P\&S\/} {\bf 3}, 131--150.

\smallskip
\nin\refs 
E.\ Bolthausen (1984) 
   An estimate of the remainder in a combinatorial central limit theorem. 
    {\it Z.\ Wahrscheinlichkeitstheorie verw.\ Geb.\/} {\bf 66}, 379--386.
 
\smallskip
\nin\refs
T.\ C.\ Brown \& A.\ Xia (2001)
 Stein's method and birth-death processes. 
 {\it Ann.\ Probab.\/} {\bf 29},  1373--1403.

\smallskip
\nin\refs
S.\ Chatterjee (2007)
 Stein's method for concentration inequalities. 
 {\it Probab. Theory Rel.\ Fields\/} {\bf 138}, 305--321.
 
\smallskip
\nin\refs
S.\ Chatterjee (2008)
 A new method of normal approximation. 
 {\it Ann.\ Probab.\/} {\bf 36},  1584--1610. 
 
\smallskip
\nin\refs
S.\ Chatterjee (2009)  
 Fluctuations of eigenvalues and second order Poincar\'e inequalities. 
 {\it Probab.\ Theory Rel.\ Fields\/} {\bf 143},  1--40. 
  
\smallskip
\nin\refs
S.\ Chatterjee (2010)
  Spin glasses and Stein's method. 
  {\it Probab.\ Theory Rel.\ Fields\/}~(to appear) 

\smallskip
\nin\refs
L.\ H.\ Y.\ Chen (1975)
  Poisson approximation for dependent trials.
    {\it Ann.\ Probab.\/} {\bf 3}, 534--545.
    
\smallskip
\nin\refs
L.\ H.\ Y.\ Chen (1986) The rate of convergence in a central limit
     theorem for dependent random variables with arbitrary index set.
     IMA Preprint Series \#243, Univ.\ Minnesota.

\smallskip
\nin\refs
L.\ H.\ Y.\ Chen \& Q.-M.\ Shao (2001)
A non-uniform Berry--Esseen bound via Stein's method.
{\it Prob.\ Theory Rel.\ Fields\/} {\bf 120}, 236--254.

\smallskip
\nin\refs
L.\ H.\ Y.\ Chen \& Q.-M.\ Shao (2004)
  Normal approximation under local dependence.
  {\it Ann.\ Probab.\/} {\bf 32}, 1985--2028.

\smallskip
\nin\refs
L.\ H.\ Y.\ Chen \& Q.-M.\ Shao (2007)
Normal approximation for nonlinear statistics using a concentration ineqality approach.
{\it Bernoulli\/} {\bf 13}, 581--599.

\smallskip
\nin\refs
P.\ Diaconis (2004)
  Stein's method for Markov chains: first examples.
  In: {\it Stein's method: expository lectures and applications,\/}
  Eds P.\ Diaconis \& S.\ Holmes, IMS Lecture Notes {\bf 46}, pp.~27--44.

\smallskip
\nin\refs
P.\ Diaconis \& S.\ Holmes (Eds.) (2004)
 {\it Stein's method: expository lectures and applications.\/}
 IMS Lecture Notes Vol.\ 46, Beachwood, Ohio.

\smallskip
\nin \refs
R.\ V.\ Erickson (1974) 
  $L_1$ bounds for asymptotic normality of $m$--dependent sums using Stein's technique.
       {\it Ann.\ Probab.\/} {\bf 2}, 522--529.

\smallskip
\nin \refs
J.\ Fulman (2004)
  Stein's method and non--reversible Markov chains.
  In: {\it Stein's method: expository lectures and applications,\/}
  Eds P.\ Diaconis \& S.\ Holmes, IMS Lecture Notes {\bf 46}, pp.~69--78.

\smallskip
\nin \refs
J.\ Fulman (2005)
 Stein's method and descents after riffle shuffles. 
 {\it Electron.\ J.\ Probab.\/} {\bf 10},  901--924.

\smallskip
\nin\refs 
J.\ Fulman (2009)
 Stein's method and characters of compact Lie groups. 
 {\it Comm.\ Math.\ Phys.\/} {\bf 288}, 1181--1201.

\smallskip
\nin\refs 
J.\ Glaz, J.\ Naus, M.\ Roos \& S.\ Wallenstein (1994)
Poisson approximations for distribution and moments of ordered $m$-spacings.
{\it J.\ Appl.\ Probab.\/} {\bf 31A}, 271--281.

\smallskip
\nin\refs 
A.\ P.\ Godbole (1993) 
  Approximate reliabilities of $m$--consecutive--$k$--out--of--$n$: Failure systems. 
  {\it Statist.\ Sinica\/} {\bf 3}, 321--328.

\smallskip
\nin \refs
F.\ G\"otze (1991) 
  On the rate of convergence in the multivariate CLT. 
  {\it Ann.\ Probab.\/} {\bf 19}, 724--739.

\smallskip
\nin\refs 
F.\ G\"otze \& A.\ N.\ Tikhomirov (2005)
Limit theorems for spectra of random matrices with martingale structure.
In: {\it Stein's method and applications,\/} Eds.\ A.\ D.\ Barbour \& L.\ H.\ Y.\ Chen,
IMS Lecture Note Series Volume 5, World Scientific Press, Singapore, pp.~181--194.

\smallskip
\nin \refs
L.\ Goldstein \& Y.\ Rinott (1996)
Multivariate normal approximations by Stein's method and size bias couplings.
{\it J.\ Appl.\ Probab.\/} {\bf 33}, 1--17.

\smallskip
\nin\refs 
L.\ Goldstein \& G.\ Reinert (1997)
  Stein's method and the zero--bias transformation with application to
  simple random sampling. {\it Ann.\ Appl.\ Probab.\/} {\bf 7}, 935--952.

\smallskip
\nin\refs
L.\ Goldstein \& G.\ Reinert (2005)
Zero biasing in one and higher dimensions, and applications.
In: {\it Stein's method and applications,\/} Eds.\ A.\ D.\ Barbour \& L.\ H.\ Y.\ Chen,
IMS Lecture Note Series Volume 5, World Scientific Press, Singapore, pp.~1--18.

\smallskip
\nin\refs 
P.\ Hall and A.\ D.\ Barbour (1984) 
   On the rate of Poisson convergence.
   {\it Math.\ Proc.\ Cam.\ Phil.\ Soc.\/} {\bf 95}, 473--480.

\smallskip
\nin \refs
S.\ T.\ Ho \& L.\ H.\ Y.\ Chen (1978)
An $L_p$ bound for the remainder in a combinatorial central limit theorem.
{\it Ann.\ Probab.\/} {\bf 6}, 231--249.

\smallskip
\nin \refs
S.\ Holmes \& G.\ Reinert (2004)
  Stein's method for the bootstrap.
  In: {\it Stein's method: expository lectures and applications,\/}
  Eds P.\ Diaconis \& S.\ Holmes, IMS Lecture Notes {\bf 46}, pp.~95--133.

\smallskip
\nin \refs
W.--L.\ Loh (1992) 
    Stein's method and multinomial approximation.
     {\it Ann.\ Appl.\ Probab.\/} {\bf 2}, 536--554.

\smallskip
\nin \refs
H.\ M.\ Luk (1994) 
  {\it Stein's method for the gamma distribution and 
     related statistical applications.\/} 
   Ph.\ D.\ Thesis, Univ.\ Southern California.

\smallskip
\nin\refs
E.\ Meckes (2008) 
  Linear functions on the classical matrix groups. 
  {\it Trans.\ Amer.\ Math.\ Soc.\/} {\bf 360}, 5355--5366. 

\smallskip
\nin\refs
E.\ Meckes (2009)
 On the approximate normality of eigenfunctions of the Laplacian. 
 {\it Trans.\ Amer.\ Math.\ Soc.\/} {\bf 361},  5377--5399.
 
\smallskip
\nin\refs
C.\ Neuhauser (1994) 
  A Poisson approximation theorem for sequence
     comparisons with insertions and deletions.  
  {\it Ann.\ Statist.\/} {\bf 22}, 1603--1629.
 
\smallskip
\nin\refs
I.\ Nourdin \& G.\ Peccati (2009a)
 Stein's method on Wiener chaos.
 {\it Probab.\ Theory Rel.\ Fields\/} {\bf 145}, 75--118. 

\smallskip
\nin\refs
I.\ Nourdin \& G.\ Peccati (2009b)
 Stein's method and exact Berry-Esseen asymptotics for 
 functionals of Gaussian fields. 
  {\it Ann.\ Probab.\/} {\bf 37},  2231--2261. 
 
\smallskip
\nin \refs
E.\ Pek\"oz (1996) 
  Stein's method for geometric approximation.
     {\it J.\ Appl.\ Probab.\/} {\bf 33},  707--713. 

\smallskip
\nin\refs
M.\ D.\ Penrose \& J.\ E.\ Yukich (2005)
Normal approximation in geometric probability.
In: {\it Stein's method and applications,\/} Eds.\ A.\ D.\ Barbour \& L.\ H.\ Y.\ Chen,
IMS Lecture Note Series Volume 5, World Scientific Press, Singapore, pp.~37--58.

\smallskip
\nin\refs
Ju.\ V.\ Prohorov (1953) 
  Asymptotic behaviour of the binomial distribution.
      {\it Uspekhi Math.\ Nauk} {\bf 83}, 135--43.

\smallskip
\nin \refs
G.\ Reinert (1995) 
  A weak law of large numbers for empirical measures via Stein's method. 
  {\it Ann.\ Probab.\/} {\bf 23}, 334--354.

\smallskip\nin\refs
G.\ Reinert \& A.\ R\"ollin (2009) 
  Multivariate normal approximation with Stein's 
  method of exchangeable pairs under a general linearity condition. 
  {\it Ann.\ Probab.\/} {\bf 37}, 2150--2173. 

\smallskip
\nin\refs
Y.\ Rinott \& V.\ Rotar (1996)
  A multivariate CLT for local dependence 
  with $n^{-1/2}\log n$ rate and applications to multivariate graph related statistics. 
  {\it J.\ Multivariate Anal.\/} {\bf 56},  333--350.

\smallskip
\nin\refs
Y.\ Rinott \& V.\ Rotar (1997)
  On coupling constructions and rates 
  in the CLT for dependent summands with applications to the antivoter model and weighted 
  $U$-statistics. 
  {\it Ann.\ Appl.\ Probab.\/} {\bf 7}, 1080--1105.

\smallskip\nin\refs
A.\ R\"ollin (2008)
  A note on the exchangeability condition in Stein's method.
  {\it Statist.\ Probab.\ Lett.\/} {\bf 78}, 1800--1806.

\smallskip
\nin \refs
M.\ Roos (1994) 
  Stein's method for compound Poisson approximation: the local approach. 
  {\it Ann.\ Appl.\ Probab.\/} {\bf 4}, 1177--1187.

\smallskip
\nin\refs
R.\ L.\ Smith (1988)
Extreme value theory for dependent sequences via the Stein--Chen method of
Poisson approximation.  {\it Stoch.\ Procs.\ Applics\/} {\bf 30}, 317--327.

\smallskip
\nin\refs
C.\ Stein (1972) 
  A bound for the error in the normal approximation to the distribution of a sum of 
  dependent random variables.
  {\it Proc.\ Sixth Berkeley Symp.\ Math.\ Statist.\ Probab.\/} {\bf 2}, 583--602.

\smallskip
\nin\refs
C.\ Stein (1986) {\it Approximate computation of expectations.\/}
     IMS Lecture Notes Vol.\ 7, Hayward, Calif.

\smallskip
\nin \refs
C.\ Stein (1992) 
  A way of using auxiliary randomization.
  In: {\it Probability Theory\/}, Eds L.\ H.\ Y.\ Chen,
     K.\ P.\ Choi, K.\ Hu and J.--H.\ Lou, W.\ de Gruyter, Berlin, pp.~159--180. 

\smallskip
\nin \refs
C.\ Stein, P.\ Diaconis, S.\ Holmes \& G.\ Reinert (2004)
  Use of exchangeable pairs in the analysis of simulations.
  In: {\it Stein's method: expository lectures and applications,\/}
  Eds P.\ Diaconis \& S.\ Holmes, IMS Lecture Notes {\bf 46}, pp.~1--26.

\smallskip
\nin\refs
H.\ Takahata (1983) 
  On the rates in the central limit theorem for weakly dependent random fields. 
  {\it Z.\ Wahrscheinlichkeitstheorie verw. Geb.\/} {\bf 64}, 445--456.

\bye